\documentclass[12pt]{article}
\usepackage{amsmath, amsthm,
amssymb,extsizes}
\DeclareMathOperator{\diag}{diag}

\title{Classification
of squared normal
operators on unitary
and Euclidean
spaces\thanks{Partially
supported by grants of
CNPq (307812/2004-9)
and FAPESP
(2005/60337-2 and
05/59407-6).} }
\author{Vyacheslav Futorny%
\\Department of Mathematics, University of
S\~{a}o Paulo\\
S\~{a}o Paulo, Brazil,
futorny@ime.usp.br
 \and
Roger A. Horn%
\\Department of Mathematics, University of
Utah\\ Salt Lake City,
Utah 84112-0090, USA,
rhorn@math.utah.edu
\and Vladimir V.
Sergeichuk
\\ Institute of
Mathematics,
Tereshchenkivska 3\\
Kiev, Ukraine,
sergeich@imath.kiev.ua}

\date{}

\begin{document}

 \maketitle

\renewcommand{\le}{\leqslant}
\renewcommand{\ge}{\geqslant}

\newtheorem{theorem}{Theorem}
\newtheorem{lemma}{Lemma}[section]

UDC 512.643 \medskip

\begin{abstract}
We give a canonical
form for a complex
matrix, whose square is
normal, under
transformations of
unitary similarity as
well as a canonical
form for a real matrix,
whose square is normal,
under transformations
of orthogonal
similarity.

\end{abstract}

For a complex matrix
$A$ such that $A^2$ is
normal, we have
previously used a
regularization
algorithm and the
theory of cosquares to
construct two
canonical forms under
transformations of
unitary similarity
$A\mapsto S^{-1}AS$
($S$ is a unitary
complex matrix)
\cite{h-s_unit}. We
now obtain these
canonical forms by
employing ideas of
Littlewood's algorithm
\cite{litt,ser_unit}.
For a real matrix $B$
whose square is
normal, we also give a
canonical form under
transformations of
orthogonal similarity
$B\mapsto R^{-1}BR$
($R$ is an orthogonal
real matrix).

Our results provide
canonical matrices of
linear operators
${\cal A}\colon U\to
U$ such that ${\cal
A}^2$ is a normal
operator and $U$ is a
unitary or Euclidean
space since changes of
the basis transform
the matrix of $\cal A$
by unitary or,
respectively,
orthogonal similarity.

The problem of
classifying linear
operators ${\cal
A}\colon U\to U$ such
that ${\cal A}^3=0$
and $U$ is a unitary
space contains the
problem of classifying
arbitrary linear
operators on unitary
spaces (see
\cite[p.\,45]{ser_unit}
or \cite{h-s_unit}).
Thus, the condition
``${\cal A}^3$ is
normal'' and even the
condition ``${\cal
A}^3=0$'' do not
simplify the problem
of classifying linear
operators ${\cal A}$
on a unitary space.

\section{Squared
normal complex
matrices under unitary
similarity}\label{sec1}

\begin{theorem}[\cite{h-s_unit}]
\label{tt1} Let $A$ be
a square complex
matrix such that $A^2$
is normal. Then $A$ is
unitarily similar to
\begin{itemize}
  \item[\rm(a)]
a direct sum of
blocks, each of which
is
\begin{equation}\label{s1}
\begin{bmatrix} \lambda
\end{bmatrix}\ \text{ or
}\
\begin{bmatrix}
\mu & r\\ 0&-\mu
\end{bmatrix},\qquad
\begin{matrix}
\lambda,\mu\in\mathbb
C,\
0\le \arg(\mu)<\pi,\\
r\in\mathbb R,\ r>0;
\end{matrix}
\end{equation}
and

  \item[\rm(b)]
to a direct sum of
blocks, each of which
is
\begin{equation}\label{s2}
\begin{bmatrix} \lambda
\end{bmatrix}\ \text{ or
}\
\tau\!\begin{bmatrix}
0 & 1\\ \nu&0
\end{bmatrix}
,\qquad
\begin{matrix}
\lambda, \nu \in
\mathbb C,\ |\nu|<1,\\
\tau\in \mathbb R,\
\tau>0.
\end{matrix}
\end{equation}
\end{itemize}
These direct sums are
uniquely determined by
$A$, up to permutation
of summands.
Conversely, if $A$ is
unitarily similar to a
direct sum of blocks
of the form \eqref{s1}
or \eqref{s2}, then
$A^2$ is normal.
\end{theorem}

\begin{proof}
(a) Let $A$ be a
squared normal complex
matrix. Let
$\lambda_1,\dots,\lambda_t$
be all the distinct
eigenvalues of $A$
ordered such that
\begin{equation}\label{kjt}
0\le
\arg(\lambda_i)<\pi\
\text{ and }\
-\lambda_i\in\{\lambda_1
,\dots,\lambda_t\}\quad\Longrightarrow
\quad
\lambda_{i+1}=-\lambda_i
\end{equation}
for every nonzero
eigenvalue
$\lambda_i$. Schur's
unitary
triangularization
theorem \cite[Theorem
2.3.1]{hor} ensures
that $A$ is unitarily
similar to a matrix of
the form
\[
T=\begin{bmatrix}
\Lambda_{\lambda_1}&T_{12}&\dots
&T_{1t}\\
&\Lambda_{\lambda
_2}&\ddots&\vdots\\
&&\ddots&T_{t-1,t} \\
0&&&\Lambda_{\lambda_
t}
\end{bmatrix},
\]
in which every
$\Lambda_{\lambda_i}$
is an $n_i\times n_i$
upper triangular
matrix of the form
\[
\begin{bmatrix}
\lambda_i && *\\
&\ddots\\
0&&\lambda_i
\end{bmatrix}.
\]

If $0$ is an
eigenvalue of $A$ and
the corresponding
block
$\Lambda_0=[a_{ij}]$
is nonzero, then the
sets of indices
$\{i\,|\,a_{ij}\ne
0\}$ and
$\{j\,|\,a_{ij}\ne
0\}$ are disjoint
because $\Lambda_0$ is
upper triangular and
$\Lambda_0 ^2=0$. We
reduce $\Lambda_0$  by
permutation similarity
transformations to the
form
\[
\Lambda_0=\begin{bmatrix}
0 & *\\0& 0
\end{bmatrix}
\]
with square diagonal
blocks.

Since $A^2$ is normal,
$T^2$ is normal, too.
But $T^2$ is upper
triangular, hence
\begin{equation*}\label{htc}
T^2=\lambda_1^2
I_{n_1} \oplus
\lambda_2^2
I_{n_2}\oplus\dots
\oplus\lambda_t^2
I_{n_t}.
\end{equation*}
This implies that
$\Lambda_{\lambda_i}
=\lambda_iI_{n_i}$ if
${\lambda_i}\ne 0.$

If $T_{ij}\ne 0$, then
$\lambda_i^2=\lambda_j^2$
since $T$ commutes
with $T^2$. By \eqref{kjt} we have
$j=i+1$, $0\le
\arg(\lambda_i)<\pi$,
and
$\lambda_{i+1}=-\lambda_i$.

Thus, $T$ is a direct
sum of matrices of two
types:
\begin{equation*}\label{kju}
\Lambda_{\lambda}=\lambda
I\quad\text{and}\quad
T_{\mu
}=\begin{bmatrix} \mu
I & F_{\mu }\\0& -\mu
I\end{bmatrix}\ (0\le
\arg(\mu)<\pi,\ F_{\mu
}\ne 0).
\end{equation*}
If $F_{\mu
}=U\Sigma_{\mu } V$ is
a singular value
decomposition, then
$T_{\mu}$ is unitarily
similar to
\[
S^{-1}T_{\mu}S=
\begin{bmatrix} \mu I
& \Sigma_{\mu }
\\0& -\mu
I\end{bmatrix},\qquad
S:=\begin{bmatrix} U &
0\\0&
V^{-1}\end{bmatrix}.
\]

Therefore, $A$ is
unitarily similar to a
direct sum of matrices
of the form
\eqref{s1}. Let us
prove that this sum is
uniquely determined by
$A$, up to permutation
of summands. We give a
direct proof, though
we could use Theorem
3.1 in
\cite{ser_unit}, which
states that every
system of linear
mappings of unitary
spaces uniquely
decomposes into a
direct sum of
indecomposable
systems, up to
isomorphism of
summands.

Let $T$ and $T'$ be
two matrices that are
direct sums of blocks
of the form
\eqref{s1}. Let
$\lambda $ be an
eigenvalue of $T$.
Grouping together the
summands with the
eigenvalues $\lambda $
and $-\lambda $, we
obtain
\[ T=T_1\oplus
T_2,\qquad
T'=T'_1\oplus T'_2,\]
in which $T_1$ and
$T'_1$ are direct sums
of blocks whose
eigenvalues are
$\lambda$ or
$-\lambda$; $T_2$ and
$T'_2$ are direct sums
of blocks that have no
eigenvalues $\lambda$
and $-\lambda$.
Suppose $T$ and $T'$
are unitarily similar.
We claim that their
decompositions into
direct sums coincide
up to permutation of
summands. It suffices
to verify that the
decompositions of
$T_1$ and $T_1'$
coincide up to
permutation of
summands. Let $S$ be a
unitary matrix such
that $TS=ST'$. Since
$T_1$ and $T_1'$ have
no common eigenvalues
with $T_2$ and $T_2'$,
the matrix $S$ has the
form $S=S_1\oplus S_2$
and
\begin{equation}\label{log}
T_1S_1=S_1T_1'.
\end{equation}
Thus, $T_1$ is
unitarily similar to
$T'_1$.

If $T_1$ is a direct
sum of $1\times 1$
blocks $[\lambda]$ and
$[-\lambda]$, then the
decompositions of
$T_1$ and $T_1'$
coincide up to
permutation of
summands. Suppose
$T_1$ has a $2$-by-$2$
direct summand of the
form \eqref{s1}.
Reduce $T_1$ and
$T_1'$ by permutation
similarity
transformations to the
form
\begin{equation*}\label{jtc}
T_1=\begin{bmatrix}
\lambda I_p & \Sigma
\\0& -\lambda
I_q\end{bmatrix},\qquad
T_1'=\begin{bmatrix}
\lambda I_{p'} &
\Sigma'
\\0& -\lambda
I_{q'}\end{bmatrix},
\end{equation*}
in which
\[
\Sigma=\diag(r_1,
\dots,r_k)\oplus
0,\qquad
\Sigma'=\diag(r'_1,
\dots,r'_{k'})\oplus 0
\]
are real matrices such
that
\[
r_1\ge\dots\ge
r_k>0,\qquad
r'_1\ge\dots\ge
r_{k'}>0,
\]
and if $\lambda =0$
then $\Sigma$ and
$\Sigma'$ have no zero
columns (hence their
columns are linearly
independent; we gather the zero columns and rows in the first vertical strip and the first horizontal strip of $T_1$ and $T'_1$). Then
\eqref{log} implies
that $S_1=U\oplus V$,
in which $U$ is
$p\times p$. Thus,
$\Sigma V=U\Sigma'$
and the uniqueness of
a singular value
decomposition ensures
that $\Sigma=\Sigma'$.

We have proved that
all direct sums of
matrices of the form
\eqref{s1} are
canonical under
unitary similarity.

(b) Let us prove that
all direct sums of
matrices of the form
\eqref{s2} are
canonical, too. It
suffices to verify
that the mapping
\begin{equation*}\label{kdt}
f\colon\quad
\tau\!\begin{bmatrix}
0 & 1\\ \nu&0
\end{bmatrix}\
\longmapsto\
\tau\!\begin{bmatrix}
\sqrt{\nu} & 1-|\nu|\\
0&-\sqrt{\nu}
\end{bmatrix}
\end{equation*}
(in which
$0\le\arg(\sqrt{\nu})<\pi$;
that is, $\sqrt{\nu}$
is the principal
square root of $\nu$)
is a bijection of the
set of matrices of the
form
\begin{equation}\label{crg}
M_{\nu,\tau}
:=\tau\!\begin{bmatrix}
0 & 1\\ \nu&0
\end{bmatrix},\qquad
\begin{matrix}
\nu \in
\mathbb C,\ |\nu|<1,\\
\tau\in \mathbb R,\
\tau>0,
\end{matrix}
\end{equation}
onto the set of
matrices of the form
\begin{equation*}\label{crg1}
N_{\mu
,r}:=\begin{bmatrix}
\mu & r\\ 0&-\mu
\end{bmatrix},\qquad
\begin{matrix}
\mu\in\mathbb
C,\ 0\le \arg(\mu)<\pi,\\
r\in\mathbb R,\ r>0,
\end{matrix}
\end{equation*}
and each
$M_{\nu,\tau}$ is
unitarily similar to
$f(M_{\nu,\tau})$.

First we prove that
$f$ is a bijection.
Fix $N_{\mu ,r}$ and
verify that it has
exactly one preimage
$M_{\nu,\tau}$. The
equality
$f(M_{\nu,\tau})=N_{\mu
,r}$ is valid if and
only if
\begin{equation}\label{jju}
\tau\sqrt{\nu}=\mu,\qquad
\tau(1-|\nu|)=r
\end{equation}
if and only if
\begin{equation}\label{hfc}
\frac{\mu}{r}
=\frac{\sqrt{\nu}}
{1-|\nu|}, \qquad
\tau=\frac{r}{1-|\nu|}.
\end{equation}
Write the complex
numbers $\mu/r$ and
$\nu$ in polar form:
$\mu/r=\rho
e^{i\varphi}$ and
$\nu=\chi e^{i\psi}$.
The first equality in
\eqref{hfc} is valid
if and only if
\begin{equation}\label{mdt}
\rho=\frac{\sqrt{\chi}}
{1-\chi}, \qquad
\varphi=\frac{\psi}2.
\end{equation}
There is a unique real
$\chi$ in the interval
$(0,1)$ that satisfies
the first equality in
\eqref{mdt} because
the real function \[y=
\frac{\sqrt{x}}
{1-x}\] steadily
increases from $0$ to
$+\infty$ on the
interval $[0,1)$.
Thus, there is a
unique $\nu$ that
satisfies the first
equality in
\eqref{hfc}. We find
$\tau $ from the
second equality and
obtain the required
preimage
$M_{\nu,\tau}$ of
$N_{\mu ,r}$.

Each matrix pair
$(M_{\nu,\tau}, N_{\mu
,r})$ with $N_{\mu ,r}
=f(M_{\nu,\tau})$ is
completely determined
by the parameters $\mu \in\mathbb C$
and $\tau \in \mathbb
R$. These parameters satisfy $0\le
\arg(\mu)<\pi$ and
$\tau >0$ since by
\eqref{jju} they
determine the
remaining parameters
$\nu$ and $r$:
\[
\nu=\frac{\mu^2}{\tau^2},
\qquad r
=\tau(1-|\nu|)=\tau-
\frac{\mu\bar\mu}{\tau}.
\]
Using these
equalities, we obtain
\[
\mu(\tau-r)=
\frac{\mu^2\bar{\mu}}{\tau}
=\nu\bar{\mu}\tau
\]
and so
\[
N_{\mu
,r}S=SM_{\nu,\tau},\qquad
S:=\frac{1}{\sqrt{\tau^2+
\mu\bar\mu}}
\begin{bmatrix}
\tau & \bar\mu\\
-\mu&\tau
\end{bmatrix}.
\]
Therefore,
$M_{\nu,\tau}$ is
unitarily similar to
$N_{\mu ,r}$.
\end{proof}

Littlewood's algorithm
\cite{litt,shap,ser_unit}
transforms each square
complex matrix $A$
into a matrix $A_{\rm
can}$ that is
unitarily similar to
$A$. Two square
matrices $A$ and $B$
are unitarily similar
if and only if they
are transformed into
the same matrix
$A_{\rm can}=B_{\rm
can}$. Thus, $A_{\rm
can}$ is a canonical
form of $A$ under
unitary similarity.
The structure of
$A_{\rm can}$ is
studied in
\cite{ser_funk,ser_unit}.
If $A$ is squared
normal, then $A_{\rm
can}$ is
permutationally
similar to a direct
sum of blocks of the
form \eqref{s1}.

\section{Squared normal real
matrices under
orthogonal
similarity}\label{sec2}

The {\it
realification} of an
$m\times n$ complex
matrix $M$  is the
$2m\times 2n$ real
matrix $M^{\mathbb R}$
obtained by replacing
every entry $a+bi$ of
$M$ by the $2\times 2$
block
\begin{equation*}\label{1j}
\begin{bmatrix}
 a&-b\\b&a
\end{bmatrix}
\end{equation*}

The real Jordan form
of $A\in\mathbb
R^{n\times n}$ can be
obtained from the
canonical Jordan form
of $A$ by replacing
all pairs of complex
conjugate Jordan
blocks
\[
J_n(a+bi)\oplus
J_n(a-bi),\qquad b>0,
\]
by $J_n(a+bi)^{\mathbb
R}$, see \cite[Theorem
3.4.5]{hor}. A real
canonical form of a
normal matrix
$A\in\mathbb
R^{n\times n}$ under
similarity can be
obtained from its
diagonal canonical
form by replacing all
pairs of complex
conjugate diagonal
entries
\[[a+bi]\oplus
[a-bi],\qquad b>0,\]
by $[a+bi]^{\mathbb
R}$, see \cite[Theorem
2.5.8]{hor}.  In the
following theorem we
show that a real
canonical form of a
squared normal matrix
$A\in\mathbb
R^{n\times n}$ under
orthogonal similarity
can be obtained in the
same way from the
canonical form of $A$
given in Theorem
\ref{tt1}(b) (and from
any other canonical
form of $A$ under
unitary similarity).

\begin{theorem}
\label{tt2} Let $A$ be
a square real matrix
such that $A^2$ is
normal. Then $A$ is
orthogonally similar
to a direct sum of
real blocks of the
form
\begin{equation}\label{ca3}
\begin{bmatrix} \lambda
\end{bmatrix},\qquad
\tau \begin{bmatrix} 0
& 1\\ \nu&0
\end{bmatrix}\ (|\nu|< 1,
\ \tau >0) ,
\end{equation}
and
\begin{equation}\label{ca4}
\begin{bmatrix}
 a & -b \\
  b&a
 \end{bmatrix}\ (b>0),
 \quad
\tau\! \begin{bmatrix}
0&0&1&0\\ 0&0&0&1\\
c&-d&0&0 \\
  d&c&0&0
\end{bmatrix}
\ (d> 0,\ c^2+d^2<1,\
\tau
>0).
\end{equation}
This direct sum is
uniquely determined by
$A$, up to permutation
of summands. It can be
obtained from the
canonical form of $A$
under unitary
similarity that is a
direct sum of matrices
of the form \eqref{s2}
by replacing all pairs
of summands
\begin{equation*}\label{ytb}
[a+bi]\oplus [a-bi]\
(b>0),\quad
\tau\!\begin{bmatrix}
0 & 1\\ c+di&0
\end{bmatrix}
\oplus\tau\!
\begin{bmatrix} 0 &
1\\ c-di&0
\end{bmatrix}\ (d>0)
\end{equation*}
by the corresponding
matrices \eqref{ca4}.
Conversely, if $A$ is
unitarily similar to a
direct sum of blocks
of the form
\eqref{ca3} and
\eqref{ca4}, then
$A^2$ is normal.
\end{theorem}

\begin{proof}
Let us prove that a
complex matrix $M$ is
squared normal if and
only if its
realification
$M^{\mathbb R}$ is
squared normal. If $M$
is represented in the
form $M=A+Bi$ with $A$
and $B$ over $\mathbb
R$, then its
realification
$M^{\mathbb R}$ is
permutationally
similar to
\[
M_{\mathbb
R}:=\begin{bmatrix}
  A&-B\\B&A
\end{bmatrix}.
\]
Since
\[
\begin{bmatrix}
  A+Bi&0\\0&A-Bi
\end{bmatrix}
\begin{bmatrix}
  I&iI\\I&-iI
\end{bmatrix}=
\begin{bmatrix}
  I&iI\\I&-iI
\end{bmatrix}
\begin{bmatrix}
  A&-B\\B&A
\end{bmatrix},
\]
we have
\begin{equation*}\label{luf}
M_{\mathbb
R}=S^{-1}(M\oplus \bar
M)S=S^{*}(M\oplus \bar
M)S,
\end{equation*}
in which
\[
S:=\frac{1}{\sqrt{2}}
\begin{bmatrix}
  I&iI\\I&-iI
\end{bmatrix}
\]
is unitary.

Thus, $M^{\mathbb R}$
is squared normal if
and only if
$M_{\mathbb R}$ is
squared normal (i.e.,
$(M_{\mathbb R}^2)^*
{M_{\mathbb R}^2}=
{M_{\mathbb R}^2}
(M_{\mathbb R}^2)^*$)
if and only if
\[
(M^2)^* {M^2}= {M^2}
(M^2)^*
\quad\text{and}\quad
({\bar M}^2)^* {\bar
M^2}= {\bar M^2} (\bar
M^2)^*
\]
if and only if $M$ is
squared normal.

A canonical form of
real matrices under
orthogonal similarity
can be obtained from
any canonical form of
complex matrices under
unitary similarity as
follows. We say that a
square complex matrix
$A$ is
\emph{decomposable}
if it is unitarily
similar to a direct
sum of square matrices
of smaller size. Let
$\cal S$ be any set of
indecomposable
canonical complex
matrices under unitary
similarity (for
example, the set of
indecomposable
matrices on which
Littlewood's algorithm
acts identically).
Each matrix $M\in\cal
S$ that is unitarily
similar to a real
matrix $R$, we replace
by $R$. Each pair
$\{M,N\}\subset\cal S$
in which $M$ is not
unitarily similar to a
real matrix and $N$ is
unitarily similar to
the complex conjugate
matrix $\bar M$, we
replace by $M^{\mathbb
R}$ or by $N^{\mathbb
R}$ ($N=M$ if $M$ is
unitarily similar to
$\bar M$). Denote the
set obtained by ${\cal
S}_{\mathbb R}$.
Theorem 4.1 in
\cite{ser_unit} about
systems of linear
mappings on unitary
and Euclidean spaces
ensures that each real
matrix $A$ is
orthogonally similar
to a direct sum of
matrices from ${\cal
S}_{\mathbb R}$ and
that this sum is
determined by $A$
uniquely up to
permutation of
summands.

Let us remove from the
set ${\cal S}$ all
matrices that are not
squared normal and
construct the set
${\cal S}_{\mathbb R}$
as above. Then ${\cal
S}_{\mathbb R}$
consists of squared
normal real matrices
and each squared
normal real matrix is
orthogonally similar
to a direct sum of
matrices from ${\cal
S}_{\mathbb R}$, which
is determined uniquely
up to permutation of
summands.

Let us prove that if
$\cal S$ is the set of
matrices \eqref{s2},
then ${\cal
S}_{\mathbb R}$ is the
set of matrices
\eqref{ca3} and
\eqref{ca4}. The real
matrices of the form
\eqref{s2} give
\eqref{ca3}. Each pair
$\{[a+bi],\,[a-bi]\}
\subset \cal S$ with
$b>0$ gives
$[a+bi]^{\mathbb R}$,
which is the first
matrix in \eqref{ca3}.
Each pair
$\{M_{c+di,\tau},\,
M_{c-di,\tau}\}
\subset \cal S$ of matrices of the form
\eqref{crg}  with
$d>0$ and $c^2+d^2<1$
gives
$M_{c+di,\tau}^{\mathbb
R}$, which is the
second matrix in
\eqref{ca4}.
\end{proof}

\end{document}